\documentclass[12pt]{amsart}

\usepackage{hyperref}

\usepackage{geometry}
\usepackage{listings}
\usepackage{mathrsfs}
\usepackage[font={small},labelfont={sf,bf}]{caption}
\usepackage{color}
\usepackage{amssymb}
\usepackage[utf8]{inputenc}
\usepackage{afterpage}
\usepackage[T1]{fontenc}
\usepackage{amsfonts}
\usepackage{mathtools}
\usepackage{amsmath}
\usepackage{verbatim}
\usepackage{bm}
\usepackage{bbm}
\usepackage{enumerate}
\usepackage{amsthm}
\usepackage{xcolor}

\usepackage{array}
\newcolumntype{L}{>{\centering\arraybackslash}m{7cm}}
\usepackage{booktabs,caption}

\newtheoremstyle{named}{}{}{\itshape}{}{\bfseries}{.}{.5em}{\thmnote{#3 }#1}
\theoremstyle{named}
\newtheorem*{namedtheorem}{Theorem}

\usepackage{relsize,exscale}

\usepackage{titling}

\theoremstyle{definition}
\newtheorem{definition}{Definition}[section]
\theoremstyle{plain}
\newtheorem{theo}[definition]{Theorem}

\newtheorem{lemma}[definition]{Lemma}

\theoremstyle{remark}

\newtheorem{rem*}[definition]{}

\newcommand*\mcupinn[2]{\vcenter{\hbox{$\mathsurround=0pt
  \ifx\displaystyle#1\textstyle\else#1\fi\bigcup$}}}

\newcommand*\mcapinn[2]{\vcenter{\hbox{$\mathsurround=0pt
  \ifx\displaystyle#1\textstyle\else#1\fi\bigcap$}}}

\begin{document}

\begin{center}
    \Large
    \textbf{A compact non-convex ancient curve shortening flow}
\end{center}

\

\centerline { Jumageldi Charyyev }

\

\begin{center}
   \fontsize{10pt}{8pt}\selectfont
    \textsc{Abstract.}
    An example of a compact, non-convex, embedded ancient solution for the curve shortening flow, which is asymptotic to Yin-Yang curve, is constructed.
\end{center}

\section{Introduction}

As with other geometric flows, the analysis of the ancient
solutions of the curve shortening flow is one of the several central problems in geometric analysis. The curve shortening flow was first introduced in \cite{mullins56} and in the same paper the first examples of 
nontrivial ancient solutions were provided, including 
Grim Reaper and Yin-Yang.
The following table lists the convexity and compactness properties of previously known embedded ancient solutions of the curve shortening flow:

\begin{table}[h]
\begin{tabular}{|c|c|L|}\hline
\bf{Compact}  & \bf{Convex} & \bf{Remarks} \\\hline
Yes & Yes & All such solutions are classified in \cite{hamilton10}: shrinking circles or ancient ``paperclip'' solutions (Angenent Ovals)  \\\hline
  No  & Yes & All such solutions are classified in \cite{bourni20}: stationary line or Grim Reaper \\\hline
  No  & No & An example of such solution with an explicit parametrization was given in \cite{nakayama94}, where the solution looks like a sine curve. A family of examples, for which the solutions have finite total curvature, was given in \cite{angenent212}. The examples are constructed through particular gluing of several Grim Reaper curves. \\\hline
\end{tabular}
\end{table}

The main result of this paper completes the missing row in the above table, by giving an example of a compact, non-convex, ancient solution.  

\begin{namedtheorem}[Main] \label{Main}
There exists a compact ancient solution of the curve shortening flow asymptotic to Yin-Yang solution as $t\rightarrow -\infty$.
\end{namedtheorem}

The same result was obtained independently in \cite{angenent21}.   

\

The strategy of the proof of the theorem is indicated below.  

\

The proof of the existence of such solution is based on constructing a sequence of the curve shortening flows $\{ \alpha_t^n \}_{n=1}^{n=\infty}$, $t \in [-T_n,0]$, $T_n \rightarrow \infty$, such that the maximum of the curvature of $\alpha_t^n$ is uniformly bounded in $n$ at any given time $t \in (-\infty,0)$. Then by smooth compactness there exists a subsequence converging to an ancient solution  $ \alpha_t^\infty $, $t \in (-\infty,0]$.

\

The closed initial curves $\{\alpha_{-T_n}^n\}_{n=1}^{n=\infty}$ are constructed by ``attaching'' a chunk of a Grim Reaper to a truncated Yin-Yang curve 
further and further away from the origin.   Since the Grim Reaper is a soliton, it is reasonable to expect that the flow $\alpha^n_t$ will have a ``tip'' region which remains close to a Grim Reaper for a long time.  Likewise, since $\alpha^n_{-T_n}$ is identical to the Yin-Yang away from the tip region, it is plausible that the portion of $\alpha^n_t$ lying away from its tip region will be very close to a Yin-Yang, for a long time. The actual proof is based on a formalization of these two heuristic ideas.

\

Here is a brief description of the approach in \cite{angenent21}.  The authors begin by constructing an approximation solution -- an evolving curve violating the curve shortening flow equation by an error which satisfies an $L^1$ spacetime bound -- by gluing Grim Reapers to truncated Yin-Yangs.
The curve shortening flow with the same initial condition remains close to the approximate solution, in the sense that, the two evolving curves co-bound a region of small area.  Using a decomposition of the curve shortening flow into graphical pieces,  $L^1$ estimates, and a contradiction argument, the authors obtain uniform curvature bounds on time intervals $[-T-\frac{1}{4},-T]$ for all sufficiently large $T$.

\

When compared, the construction of older and older solutions
 and an extraction of an ancient solution as a subsequential limit is standard and the same for the two approaches.
However,
 apart from the standard ingredients such as parabolic regularity and the fact that the number of inflection points is non-increasing in time, there is a little overlap between the two arguments. Both the form of the estimates and the methods for proving them are different.  


%
%
%
%
\

{\bf Organization of the paper.}  First, in  Section \ref{Y-Y} the properties and the asymptotics of the Yin-Yang curve are introduced. Section \ref{AC} contains the construction of the initial curves.  Section \ref{HS} defines a class of barrier curves. Section \ref{MA}   formulates and proves the ``stability'' properties mentioned in the preceding paragraph  (Theorem \ref{EDL}), as well as the proof of Main Theorem. 

\

{\bf Acknowledgments.} The author wishes to thank his advisor Prof. Bruce Kleiner for his support and guidance during the course of this work.

\section{The Yin-Yang solution} \label{Y-Y}
This section presents some properties of the Yin-Yang solutions and introduces some notation which will be used later.

\

$\bf{\hat{\gamma_t}:}$ \indent Yin-Yang curve shortening flow -- a spiral-like counterclockwise rotating soliton with the unit angular speed centered at the origin and arclength parametrized on the complex plane as follows:
\begin{eqnarray} \nonumber
\hat{\gamma_t}(s) &=& e^{it}F(s) \\ \nonumber
F(s) &=& x(s)+iy(s) \\ \nonumber
x'(s) &=& \cos((x(s)^2+y(s)^2)/2) \\  
y'(s) &=& \sin((x(s)^2+y(s)^2)/2) \\ \nonumber
T(s) &=& x'(s)+iy'(s) = F'(s) \\ \nonumber
N(s) &=& -y'(s)+ix'(s) = iT(s) \\ \nonumber
H (s) &=& x(s)x'(s)+y(s)y'(s) 
\end{eqnarray}
where $s$ is clearly an arclength parametrization such that $F(0)=0$, $H(s)>0$ for $s>0$, and $N(s)$ is the inward unit normal vector. This solution first appeared in \cite{mullins56}. All the self-similar solutions of the curve shortening flow were classified in \cite{halldorsson10}.

\

One of the key properties of the curve is that $H = \langle F,T \rangle$, which justifies the fact that Yin-Yang is a rotating soliton: 
\begin{eqnarray} \nonumber
& \langle \partial_t \hat{\gamma_t}(s) , N(\hat{\gamma_t}(s)) \rangle = \langle ie^{it}F(s) , e^{it}N(s)  \rangle = \langle  iF(s) , N(s) \rangle = & \\
& = \langle F(s) , -iN(s)  \rangle = \langle F(s) , T(s)  \rangle = H(s) & \\ \nonumber
& \langle \partial_{ss} \hat{\gamma_t}(s) , N(\hat{\gamma_t}(s)) \rangle = \langle e^{it} F''(s) , e^{it} N(s) \rangle = \langle  F''(s) , N(s) \rangle = H(s) &
\end{eqnarray}
Here for $z,w \in \mathbb{C}$ the inner product is $\langle z,w \rangle := \Re(z\bar{w})$, which coincides with the usual dot product of the real vectors. The fact that $H = \langle F , T \rangle$ follows quickly from the following identities:
\begin{eqnarray} \nonumber
& x'' = -(xx'+yy') \sin ((x^2+y^2)/2) = -(xx'+yy')y' & \\
& y'' = (xx'+yy') \cos ((x^2+y^2)/2) = (xx'+yy')x' & \\ \nonumber
& H = x'y'' - x''y' = xx' + yy' = \langle F , T \rangle &
\end{eqnarray}

\

The portions of the curve with the positive and negative curvature occasionally will be referred as the $\bf{``positive \ wall"}$ and the $\bf{``negative \ wall"}$, respectively.

\

In \cite{halldorsson10} it was proved that $H(s)$ is an odd function and $H(s)>0$ for $s>0$. The function has one critical point in $[0,\infty)$, at which it has the global maximum; afterwards $H(s)$ is decreasing asymptotically to 0 and is convex in the ray $[s_i,\infty)$ for some $s_i>0$ (a positive inflection point of $H(s)$ function). Thus, $H'(s)<0$ for $s \in [s_i,\infty)$ and $H'(s) \xrightarrow {s \rightarrow \infty} 0$ .

\

From the parametrization of Yin-Yang at time $t=0$, it can be observed that $\hat{\gamma}_0$ is an integral curve of the unit vector field $V(x,y) = \left( \cos((x^2+y^2)/2),\sin((x^2+y^2)/2) \right)$ which passes through the origin. The integral curves of the vector field $V(x,y)$ foliate the plane. From the above equations it can be observed that each integral curve is another rotating soliton (counter-clockwise with unit angular speed) of the curve shortening flow. Among all these integral curves only Yin-Yang is symmetric with respect to the origin. For the different characterization of these curves see \cite{halldorsson10}.

\

$\bf{\hat{D}_t:}$ \indent the $\bf{``corridor"}$ -- the closure of one of the connected components of $\mathbb{R}^2 \setminus \hat{\gamma_t}$, for which $N(\hat{\gamma}_t(s))$, the unit normal of $\hat{\gamma_t}$ at $s$, points towards the inside of $\hat{D}_t$ when $s>0$. So, $\langle  F,N \rangle < 0$ when $s>0$ and $\langle  F,N \rangle > 0$ when $s<0$. Since $\hat{\gamma}_t$ is a proper embedding of $\mathbb{R}$ into $\mathbb{R}^2$, then by the Jordan separation theorem, the complement has two connected components; the boundary of each is the image of $\hat{\gamma}_t$ .

\

$\bf{w(s):}$ \indent defined for $s<0$, is the radial width of the corridor $\hat{D}_t$ at the point $\hat{\gamma}_t(s)$ -- the length of the line segment $[\hat{\gamma}_t(s), \hat{\gamma}_t(\sigma)]$, where $\sigma$ is such that the latter line segment lies completely inside the corridor and passes through the origin when extended. Such $\sigma$ exist uniquely because it is the first instance when the ray $[0,\hat{\gamma}_t(s))$ intersects the positive wall after $\hat{\gamma}_t(s)$.

\

\begin{lemma} \label{Asymptotics}
As $|s| \rightarrow \infty$ the mutual ratios of the following quantities 
$$ (3|s|)^{1/3}, \indent \dfrac{1}{|H(s)|}, \indent |F(s)|, \indent \dfrac{\pi}{w(s)} $$ converge to 1 i.e., the quantities are asymptotic.
\end{lemma}

\

\begin{proof}[\bf{Proof of Lemma \ref{Asymptotics}}]

It is straightforward to justify that the following identities hold true:
\begin{eqnarray} \nonumber
& x'' = -Hy', \indent y'' = Hx'  & \\ \nonumber
& (x'y-xy')' = x''y - xy'' = -H^2 & \\ \nonumber
& x'y - xy' = -\int_0^sH^2 & \\ \nonumber
& x'y - xy' = \langle F,N \rangle = -\int_0^sH^2 & \\ \nonumber
& H' = 1 + xx'' + yy'' = 1 + H(x'y - xy') & \\ 
& H' = 1 - H\int_0^sH^2 = 1 + H \langle F,N \rangle & \label{diff eqn for H}
\end{eqnarray} 

\

Consider the case when $s_i<s$. Since $H'(s) \xrightarrow {s \rightarrow \infty} 0$, from the relation (\ref{diff eqn for H}) it will follow that $1/H(s)$ and $\int_0^sH^2$ are asymptotic (i.e., there ratio converges to 1). Since both functions are increasing in $s$ and are asymptotic, then their respective derivatives, $-H'/H^2$ and $H^2$, must be asymptotic. So, $-H' \cong H^4$ (i.e., the ratio converges to 1), which implies that $H(s) \cong (3s)^{-1/3}$. The similar result is true for $s<0$.


\

Since $\partial_s(\frac{1}{2}|F(s)|^2)= \langle F(s),T(s) \rangle = H(s)$, after integrating it follows that $\frac{1}{2}|F(s)|^2 \cong \frac{1}{2}(3|s|)^{2/3}$ and hence the desired outcome $|F(s)| \cong (3|s|)^{1/3}$.

\

For the lifted polar angle $\theta = \arctan \dfrac{y}{x}$ the following identity is true:

\begin{equation}
\partial_s \theta = \dfrac{xy' - x'y}{x^2+y^2} = - \dfrac{\langle F,N \rangle}{|F|^2}
\end{equation}

\

For the given $s>0$ let $\sigma(s)>0$ be such that $\theta(\sigma(s)) - \theta(s) = 2\pi$ (i.e., one full circular turn along the positive wall). It should be noted that $\sigma(s) = s + O(s^{1/3})$, because of the almost one circular turn of radius $|F(s)|$. Then one has the following:

\begin{eqnarray} \nonumber
& 2\pi = \mathop{\mathlarger{\int}}_s^{\sigma(s)} \partial_s\theta = - \mathop{\mathlarger{\int}}_s^{\sigma(s)} \dfrac{\langle F,N \rangle}{|F|^2} 
= \mathop{\mathlarger{\int}}_s^{\sigma(s)} \dfrac{\sqrt{|F|^2-\langle F,T \rangle^2}}{|F|^2} 
= \mathop{\mathlarger{\int}}_s^{\sigma(s)} \dfrac{1}{|F|} \sqrt{1-\dfrac{H^2}{|F|^2}} &\\ & \cong \mathop{\mathlarger{\int}}_s^{\sigma(s)} \dfrac{1}{|F|} \cong \mathop{\mathlarger{\int}}_s^{\sigma(s)} (3\tau)^{-1/3} d\tau = \frac{1}{2} \left[ (3\sigma(s))^{2/3}-(3s)^{2/3} \right] = & \\ \nonumber & = \frac{1}{2} \left[ (3\sigma(s))^{1/3}+(3s)^{1/3} \right] \left[ (3\sigma(s))^{1/3}-(3s)^{1/3} \right] \cong & \\ \nonumber & \cong (3s)^{1/3} \left[ (3\sigma(s))^{1/3}-(3s)^{1/3} \right] \cong  (3s)^{1/3} |F(\sigma(s)) - F(s)| & 
\end{eqnarray} 
where $a \cong b$ means that $a$ and $b$ are asymptotic, so that $|F(\sigma(s)) - F(s)| \cong 2\pi (3s)^{-1/3}$. The approximations above are easily justifiable, because the integrands are asymptotic as $s \rightarrow \infty$. Now for the radial width it can be deduced that $w(s) \cong \pi (3s)^{-1/3}$, because $|F(\sigma(s)) - F(s)| \cong 2w(s)$. The lemma is proved.

\end{proof}

\section{Admissible Curves} \label{AC}

In this section particular examples of the curves will be constructed, which look like a portion of Grim Reaper glued to Yin-Yang. These curves will be referred to as $\bf{admissible}$ $\bf{initial}$ $\bf{conditions}$ and any curve obtained from them by the curve shortening flow will be called an $\bf{admissible}$ $\bf{curve}$.

\

Before proceeding, it should be noted that a {\bf Grim Reaper curve} 
is a curve congruent to the graph of the function $y = -\frac{1}{\lambda}\log \cos \lambda x$ for some $\lambda > 0$. The parameter $\lambda$
corresponds to the speed of the translation of the corresponding Grim Reaper solution of the curve shortening flow. The case $\lambda = 1$ is referred as the standard Grim Reaper (curve).

\

Using the notation and terminology from the previous section, the construction of admissible curves  goes as follows:

\

\begin{enumerate}

\item Take $s<0<\sigma$ such that $w(s) = |F(s)-F(\sigma)|$ (two points on the opposite walls of the corridor) and the line segment $[F(s),F(\sigma)]$ is radial. For the sake of definiteness assume that the ray $OA$ coincides with the positive $y-$axis. Denote $A=F(s)$ and $B=F(\sigma)$. $|s| \gg 1$ is chosen to be big enough, which will be specified later.

\

\item Let $a$ and $b$ be the lines perpendicular to the ray $OA$ at the points $A$ and $B$, respectively. It should be noted that, since the function $|F(\bar{s})|$ is increasing as $|\bar{s}|$ increases, there exists $\bar{\sigma}>0$ such that $F(\bar{\sigma})$ lies on the line $a$, and the portion of the positive wall between the points $F(\sigma)$ and $F(\bar{\sigma})$ lies entirely between the parallel lines $a$ and $b$. Moreover, the portion of the negative wall between the points $F(s)$ and $F(0)=O$ lies below the line $a$.

\

\item Let the curve $\rho$ be the Grim Reaper such that $(i)$ it is a graph over the line segment $AB$; $(ii)$ has the tip at the midpoint of the line segment $AB$; $(iii)$ has the width equal to $|AB|$; $(iv)$ is expanding in the direction of $-iF(s)$ i.e., in this case, in the positive $x-$axis direction.

\

\item Translate the tip of the curve $\rho$ and itself towards the origin along the line segment $AB$ until it touches the negative wall. Denote the new Grim Reaper to be $\bar{\rho}$. Let $\bar{\rho}$ be tangent to the negative wall at the point $P = F(s_0)$, $s<s_0<0$.

\

\item Let the curve $\bar{\rho}$ intersect the positive wall for the first time at the point $Q = F(\sigma_0)$, $0 < \sigma_0 < \sigma$.

\

\item Make a convex $C^2$ ``gluing'' of Grim Reaper $\bar{\rho}$ with the positive wall near the point $Q$ as follows: $(i)$ choose the point $U$ on $\bar{\rho}$ such that the distance from $U$ to $Q$ along $\bar{\rho}$ is equal to $w(s)$ and such that the point $B$ is closer to $U$ than to $Q$; $(ii)$ choose the point $V$ on the positive wall such that the distance from $V$ to $Q$ along the positive wall is also equal to $w(s)$ and now $B$ is closer to $Q$ than to $V$ (after rescaling by the factor of $\pi/w(s)$ the distances will be equal to $\pi$); $(iii)$ let the tangent line $u$ at $U$ tangent to $\bar{\rho}$ and the tangent line $v$ at $V$ tangent to the positive wall intersect at the point $W$; $(iv)$ let $\bar{\Lambda}(s)w(s)/\pi$ be equal to the distance between $Q$ and the ray $OB$; $\bar{\Lambda}(s) \rightarrow \infty$ when ${|s|\rightarrow \infty}$, as it will be explained below; $(v)$ make a $C^2$ gluing of $\bar{\rho}$ and the positive wall at the points $U$ and $V$ such that the new curve connecting $U$ and $V$ lies completely inside the triangular-like region $WUV$ bounded by the lines $u$ and $v$, and the portions of $\bar{\rho}$ and the positive wall; moreover, the gluing must be such that the curvature of the new curve connecting the points $U$ and $V$ is bounded by 1 after rescaling by $\pi/w(s)$; such gluing is possible, because after rescaling the curvatures at the points between $U$ and $V$ on $\bar{\rho}$ and the positive wall converge to zero uniformly as $|s| \rightarrow \infty$, and the angle between the tangent lines $u$ and $v$ is converging to $\pi$ as $|s| \rightarrow \infty$. This specifies how large $|s|$ should be.

\

\item Similar to the above construction, using the partition of unity make a $C^2$ gluing of $\bar{\rho}$ with the negative wall near the point $P$ such that the curvature of the newly obtained curve near the point $P$ is bounded by 1 after rescaling by $\pi/w(s)$, and such that there is only one inflection point (it can be even arranged so that the point $P$ itself is the only inflection point). This also specifies how large $|s|$ should be.

\end{enumerate}


By the virtue of the above construction, it should be noted that, the curvature of every admissible initial condition is bounded by some universal constant, say, $M$ after rescaling by the factor of $\pi/w(s)$. 

\

Since $\hat{\gamma}_0(s)$ has a unique inflection point at $s=0$, an admissible initial condition $\alpha$ also has an inflection point at the origin. Let $\{\alpha_t\}_{t \in [0,T)}$ be the evolution of the curve $\alpha = \alpha_0$ under the curve shortening flow. Clearly, $\alpha_t$ lies inside $\hat{D}_t$ due to the avoidance principle of the curve shortening flow. Because the number of the inflection points decreases under the curve shortening flow, the curves $\alpha_t$ have at most two inflection points.

\

In the above construction let $\Lambda(s)w(s)/\pi$ be equal to the distance between the point $P$ and the ray $OA$; after rescaling by the factor of $\pi/w(s)$ the quantity $\Lambda(s)$ represents the ``height'' of the portion of the attached Grim Reaper curve $\bar{\rho}$, which is closer to the negative wall (after rescaling the Grim Reaper has the standard width $\pi$). Clearly, $\Lambda(s) \rightarrow \infty$ as $|s| \rightarrow \infty$, because after rescaling the distance between the parallel lines $a$ and $b$ becomes exactly $\pi$ and the distance between the points $A$ and $P$ is now $\Lambda(s)$, the rescaled curve $\bar{\rho}$ is now the standard Grim Reaper, and in the vicinity of the point $A$ the negative wall is getting exponentially close to the circle centered at the origin of the radius $\pi|F(s)|/w(s)$ as $|s| \rightarrow \infty$. It can be further estimated that $\Lambda(s)e^{\Lambda(s)} \cong \pi^2/w(s)^2$ and hence after rescaling the curvatures near the points $P$ and $Q$ of the Grim Reaper and Yin-Yang are asymptotic to $e^{-\Lambda}$ and $e^{-\Lambda}/\Lambda$, respectively. An analogous result can be derived for $\bar{\Lambda}(s)$, the ``height'' of the portion of the attached Grim Reaper curve $\bar{\rho}$, which is closer to the positive wall i.e., the distance from the point $Q$ to the ray $OA$, where $\bar{\Lambda}(s)^2 e^{\bar{\Lambda}(s)} \cong 2 \pi^2/w(s)^2$. It can be easily deduced that $\Lambda(s) \cong \bar{\Lambda}(s)$. The estimation of the asymptotics of $\bar{\Lambda}(s)$ can be found in the Appendix.

\


Let $\hat{\gamma}^{\theta}_t := e^{i\theta}\hat{\gamma}_t$ for $t \geq 0$ and for $\theta \in (0,\pi)$, the rotated families of Yin-Yang curves. The family $\{\hat{\gamma}^{\theta}_t\}_{\theta \in (0,\pi)}$ foliates the interior of $\hat{D}_t$. For each $\theta \in (0,\pi)$ the curve $\hat{\gamma}^{\theta}_t$ intersects the curve $\alpha_t$ in at most two points, since initially $\hat{\gamma}^{\theta}_0$ intersects $\alpha_0$ in exactly two points.

\

Consider the angle function $\theta : \alpha_t \mapsto S^1$, which has a lift $\hat{\theta} : \alpha_t \mapsto \mathbb{R}$, since $\alpha_t$ is the boundary of the disc lying inside the corridor $\hat{D}_t$. The curvature function $H$ of the admissible initial condition $\alpha_0$ has exactly two zeros and both are transverse. Moreover, $\alpha_0$ has exactly two radial tangencies.

\

\begin{lemma} \label{Tip}
Let $[0,t_0)$ be the maximal time interval such that for any $t$ in that interval the curvature function $H$ of $\alpha_t$ has exactly two transverse zeros. Then for every $t \in [0,t_0)$ there are exactly two radial tangencies of $\alpha_t$, which are the non-degenerate critical points of the angle function $\hat{\theta}$ (max-min).
\end{lemma}

\

\begin{proof}[\bf{Proof of Lemma \ref{Tip}}]
Assume the contrary and let $[0,t_1)$ be the maximal time interval for which the assertion is true, where $t_1 < t_0$. A local max or min of the function $\hat{\theta}$ must be non-degenerate, because the curvature is nonzero at that point (otherwise, if the curvature vanishes, then it has a transverse zero at that point and locally there $\alpha_t$ will lie on both sides of the radial tangency, hence contradicting the local max-min assumption). From the definition of $t_1$ it follows that there must be a radial tangency at some point $p \in \alpha_{t_1}$, where $p$ is not a non-degenerate local max or min. So, $p$ is in fact an inflection point. However, then there should be a pair of radial tangencies close to $p$ on $\alpha_{t_1}$, contradicting the definition of $t_1$. Thus, $t_1 = t_0$. The lemma is proved.
\end{proof}

\

From Lemma \ref{Tip} it will follow that at time $t=t_0$ the admissible curve $\alpha_{t_0}$ has exactly one non-transverse zero of the curvature function $H$ and for further times $t>t_0$ the curves $\alpha_t$ become strictly convex.

\

$\bf{tip(\alpha_t)}:$ \indent the tip of the admissible curve $\alpha_t$, the non-degenerate critical point corresponding to the max of the angle function $\hat{\theta}$.  The uniqueness follows from the lemma when $t<t_0$, and from the convexity of $\alpha_t$ when $t\geq t_0$. 

\

$\bf{tip(\alpha_t)^{\pm}:}$ \indent the radial projections of $tip(\alpha_t)$ on the positive and negative walls of $\hat{\gamma_t}$.

\

$\bf{L_{\pm}(t):}$ \indent the distances from the origin along the walls of $\hat{\gamma}_t$ of the points $tip(\alpha_t)^{\pm}$ i.e., $tip(\alpha_t)^{+} = \hat{\gamma}_t( L_{+}(t))$ and $tip(\alpha_t)^{-} = \hat{\gamma}_t( -L_{-}(t))$.

\

$\bf{w(t)}:$ \indent the ``width'' at the tip, the radial distance between $tip(\alpha_t)^+$ and $tip(\alpha_t)^-$ or $w(-L_-(t))$ in correspondence with the notation from the previous section.

\



\section{Hyperbolic Supersolution} \label{HS}

A family of the open curves $\Gamma(s,t)$, where $t \in [0,T)$ and $s \in [-L(t), R(t)]$ for some functions $L(t), R(t) \in C^1([0,T))$, is defined to be a $\bf{hyperbolic \ supersolution}$, if it has the following properties (where $\hat{\gamma}_t$ and $\hat{D}_t$ are as in Section~\ref{Y-Y}):

\

\begin{enumerate}
\item $\Gamma(s,t) = e^{it}(F(s)+\phi(s,t)N(s))$ and the curve is completely inside the region $\hat{D}_t$ (a graph over Yin-Yang $\hat{\gamma}_t$ at each time $t$);

\

\item $\phi(s,t) = \lambda^- w(t) \exp[-a(t)(L(t)+s)] + \lambda^+ w(t) \exp[-a(t)(R(t)-s)]$, where $w(t), a(t) \in C^1([0,T))$ and $\lambda^{\pm} > 0$ (a hyperbolic function over the Yin-Yang);

\

\item $ \langle \partial_t \Gamma , \tilde{N} \rangle \geq \tilde{H} $, where $\tilde{N}$ and $\tilde{H}$ are the unit inward normal and the curvature of the $\Gamma-$curve, respectively (a supersolution to the curve shortening flow).
\end{enumerate}

\

From the definitions of the unit normal and the curvature the following statement is true:
\begin{equation}
\langle \partial_t\Gamma, \tilde{N} \rangle \geq \tilde{H} \Longleftrightarrow  \langle \partial_t\Gamma, \dfrac{i\partial_s\Gamma}{|\partial_s\Gamma|} \rangle \geq \dfrac{\langle \partial_s\partial_s\Gamma, i\partial_s\Gamma \rangle}{|\partial_s\Gamma|^3}
\end{equation}
The last inequality, the supersolution condition, can be rewritten as:
\begin{eqnarray} \label{supersolution condition}
& |\partial_s\Gamma|^2 \langle \partial_t\Gamma,i\partial_s\Gamma \rangle - \langle \partial_s\partial_s\Gamma, i\partial_s\Gamma \rangle \geq 0 & \\ \label{terms}
& A(\phi) + Q(\phi) + P(\phi) \geq 0 &
\end{eqnarray}
Where $A$, $Q$, and $P$ contain solely the linear, quadratic, and higher order terms in $\phi$, respectively.
In the Appendix it is computed that 
\begin{eqnarray}
A(\phi) &=& \phi_t - \phi_{ss} + \phi_s \langle F(s), N(s) \rangle - \phi H(s)^2 \\
Q(\phi) &=& - H(3\phi A(\phi) + 2 \phi \phi_{ss} + \phi_s^2 + \phi^2 H^2)
\end{eqnarray}

\

\section{The Main Argument} \label{MA}

This section gives the main estimate on the evolution of admissible curves, and uses it to prove the Main Theorem.

\begin{theo} \label{EDL}
For any $\epsilon > 0$ there exists $R_{\epsilon} \gg 1$ such that if $\alpha_0$ is an admissible initial condition with $\Vert tip(\alpha_0) \Vert > R_{\epsilon}$ and $\alpha_t$ is the evolution of $\alpha_0$ under the curve shortening flow with $\Vert tip(\alpha_t) \Vert \geq R_{\epsilon}$ for $t \in [0,T)$, then 

\

\begin{enumerate}

\item \label{stability} (Stability near the tip)\; In the ``neighborhood of the tip$(\alpha_t)"$ 
{i.e., the connected component of $\alpha_t$ containing $tip(\alpha_t)$ inside the ball of the radius $1000w(t)$ centered at tip$(\alpha_t)$}, the curve $\alpha_t$ is $\epsilon-$close to the standard Grim Reaper in the pointed $C^{\infty}-$topology after rescaling the curve by the factor of $\pi/w(t)$ ;

\

\item \label{proximity} (Exponentially close to the Yin-Yang) \; There exists a hyperbolic supersolution $\Gamma(s,t)= e^{it}(F(s)+\phi(s,t)N(s))$ lying inside the closed curve $\alpha_t$ for $t \in [0,T)$ and $s \in [-L(t),R(t)]$ (the ``barrier curve''), where 

\begin{eqnarray}
\phi(s,t) = \lambda^- w(t) \exp[-a(t)(L(t)+s)] + \lambda^+ w(t) \exp[-a(t)(R(t)-s)] \\ \nonumber \\ 
a(t) = \dfrac{\pi (1-\epsilon)}{w(t)}, \indent \indent L(t) = L_-(t)-\lambda w(t), \indent R(t) = L_+(t)-\lambda w(t)
\end{eqnarray}
for some $\lambda^{\pm}, \lambda > 0$ and $w(t)$ is the width at the tip.

\end{enumerate}

\
{Moreover, there exists a constant $M_{\epsilon} \gg 1$ such that the curvature of $\alpha_t$ away from the neighborhood of the tip (see the condition (\ref{stability}) above) is bounded by $M_{\epsilon}$ after rescaling by the factor of $\pi/w(t)$ as long as the conditions (\ref{stability}) and (\ref{proximity}) hold i.e., for $t \in [0,T)$ (the control of the curvature away from the tip).}

\end{theo}

\

\begin{proof}[\bf{Proof of Main Theorem \ref{Main}}]

Fix any $\epsilon>0$ and consider a sequence of time-shifted admissible curves $\{ \alpha_t^n \}_{n=1}^{n=\infty} $, $t \in [-T_n,0) $ (the extinction occurs at $ t=0 $), with $ \Vert tip(\alpha_{-T_n}^n) \Vert \rightarrow \infty $ (and hence $ T_n \rightarrow \infty $). Let $ \{ T_n^{\epsilon} \}_{n=1}^{n=\infty} $ be the sequence of times such that $ \Vert tip(\alpha_{-T_n^{\epsilon}}^n) \Vert = R_{\epsilon} $. 

\

Clearly, $ \{ T_n^{\epsilon} \}_{n=1}^{n=\infty} $ must be bounded, because the admissible curves will be sweeping an area of size $ 2 \pi T_n^{\epsilon} $ starting from the moment when they cross the threshold $ R_{\epsilon} $ until their extinction, which is clearly a finite area. Passing to a subsequence it can be taken that $T_n^{\epsilon} \xrightarrow {n \rightarrow \infty} T^{\epsilon} $.

\

For any fixed $ t \in (-\infty,-T^{\epsilon}) $ the sequence $ \{ \Vert tip(\alpha_{t}^n) \Vert \}_{n=1}^{n=\infty} $ must be bounded. Otherwise, there will be a sequence of admissible curves defined for the time interval $[t,0]$, where each curve goes extinct in finite amount of time $-t$ and sweeps an area of size $ -2\pi t $ in the process, while the tips of the curves at the time $t$ can get arbitrarily far away from the origin. Because of exponential proximity of the admissible curve $\alpha_{t}^n$ to the walls of $\hat{D_t}$, the area of the region enclosed by the hyperbolic supersolution or the barrier curve $ \Gamma(s,-t) $, after connecting its endpoints by the line segment, gets arbitrarily large as $n \rightarrow \infty$, since $ \Vert tip(\alpha_{t}^n) \Vert \rightarrow \infty $. For every $n$ this region lies completely inside $\alpha_{t}^n$, which contradicts with the above fact of $\alpha_{t}^n$ enclosing a region of the area $ -2 \pi t $ for any $n$. So, say, $ \Vert tip(\alpha_{t}^n) \Vert \leq R^t < \infty $ for any $n$.

\

Fix $ t \in (-\infty,-T^{\epsilon}) $. In the vicinity of their tips the curvatures of $\alpha_{t}^n$ are bounded by $ (1+\epsilon)R^t $ for any any $n$, because of being $\epsilon-$close to the standard Grim Reaper after rescaling by the factor of $R^t$. {Away from the tips the maximum of the curvatures of $\alpha_{t}^n$ are bounded uniformly in $n$, in fact, they are bounded by $M_{\epsilon}R^t$.}

\

Thus, for every $ t \in (-\infty,-T^{\epsilon}) $, the maximum of the curvature function of $\alpha_{t}^n$ is bounded uniformly in $n$. So, there exists a subsequence of $\{ \alpha_t^n \}_{n=1}^{n=\infty} $ which converges to some curve shortening flow $\chi_t$, $ t \in (-\infty,-T^{\epsilon}) $. The curve shortening flow $\chi_t$, $t \in (-\infty,0)$, is clearly ancient, non-convex, and asymptotic to Yin-Yang as $t \rightarrow -\infty$. The theorem is proved.

\end{proof}

\

\begin{proof}[\bf{Proof of Theorem \ref{EDL}}]
For the hyperbolic supersolution it can be further specified that one can take $\lambda^{\pm} = \frac{1}{\pi (1-\epsilon)}$ and $\lambda=\frac{100}{\pi (1-\epsilon)}$, so that the function $\phi$ will have a simpler form:
\begin{eqnarray}
\phi(s,t) = \dfrac{1}{a} \exp[-a(L+s)] + \dfrac{1}{a} \exp[-a(R-s)] \\
a = \dfrac{\pi (1-\epsilon)}{w}, \indent \indent L = L_--\dfrac{100}{a}, \indent R = L_+-\dfrac{100}{a}
\end{eqnarray}

Assume the contrary. Then for some $\epsilon>0$ there is a sequence of the admissible initial conditions $\{\alpha^j_0\}_{j=1}^{j=\infty}$ with the corresponding evolutions $\{\alpha^j_t\}_{j=1}^{j=\infty}$, $t \in [0,T_j]$, under the curve shortening flow, where the time $T_j$ is the first instance when one of two properties  {in Theorem \ref{MA}}, either ``stability'' or ``proximity'', fail and one has that $\Vert tip(\alpha_{T_j}^j) \Vert =: R_j \rightarrow \infty$. {This leads to the analysis of two cases.}

\

{\bf{Case 1$^{\circ}$:}} Assume that the stability fails for infinitely many indices $j$ (hence after passing to the subsequence can be assumed for all $j$). Consider $\beta^j_t$, the rescaled by the factor of $\pi/w(T_j)$ and time-shifted curve shortening flows of $\alpha_t^j$,  in the time intervals $[-\pi^2 T_j/w(T_j)^2,0]$. Passing to a subsequence it can be argued that $-\pi^2 T_j/w(T_j)^2 =: -K_j \rightarrow -K$ for some $0 \leq K \leq \infty$. 

\

{\bf{ \bf{Case 1$^{\circ}$a:} }} First consider the case when $K = \infty$. {Below it will be proved that,} {at any given time $t \in (-\infty,0)$ the curves $\{\beta_t^j\}_{j=1}^{j=\infty}$ have uniform bounds on their curvature: the uniform bound near the tip due to the stability combined with the uniform bound away from the tip owing to the fact of exponential proximity.} 

\

{For the sake of contradiction, assume that away from their tips the maxima of the curvature of the curves $\beta_t^j$ is an unbounded sequence in $j$. Moreover, let $-\hat{K_j}$ be the first instance when the maximum of the curvature of $\beta_t^j$ is at least $j$. There are two possibilities: either $\limsup_{j \rightarrow \infty}(K_j-\hat{K_j})>0$ (``long enough flow'') or $\limsup_{j \rightarrow \infty}(K_j-\hat{K_j})=0$ (``rapid increase of curvature''). For any $j$ at the time $-\hat{K_j}$ the part of the curve $\beta_{-\hat{K_j}}^{j}$ away from its tip is graphical over Yin-Yang, and it is confined between Yin-Yang and the corresponding hyperbolic supersolution. After interpolating between Yin-Yang and the barrier curve one gets uniform in $j$ bound on the $C^1-$norms of the curves $\beta_{-\hat{K_j}}^{j}$.}

\

{If $\limsup_{j \rightarrow \infty}(K_j-\hat{K_j})>0$, then, passing to a subsequence if necessary, the uniform in $j$ bound on $C^1-$norm combined with the parabolic regularity will yield a uniform in $j$ bounds on $C^k-$norms. This will contradict with the fact that the maxima of the curvatures of the curves $\beta_{-\hat{K_j}}^{j}$ away from the tips is unbounded in $j$.}

\

{In the case of $\limsup_{j \rightarrow \infty}(K_j-\hat{K_j}) = 0$, the uniform in $j$ bounds of $C^k-$norms of the curves $\beta_{-K_j}^j$ away from the tips and the uniform in $j$ bound of $C^0-$norm of the curves $\beta_{-\hat{K_j}}^{j}$ away from the tips will yield uniform in $j$ bounds of $C^k-$norms of the curves  $\beta_{\tau}^j$ for $\tau \in [-K_j,-\hat{K_j}]$. This causes another contradiction.} 

\

{Thus, away from their tips the maxima of the curvature of the curves $\beta_t^j$ must be a bounded sequence in $j$.}

\

Due to the uniform in $j$ bound on the curvature, the pointed sequence of the curve shortening flows $\beta_t^j$, pointed at the tips, converge in Cheeger-Gromov-Hausdorff sense to some non-compact ancient solution, say $\chi_t$. Clearly,  all the curves $\chi_t$ must be restricted inside some strip of width $\pi$. The curve $\chi_t$ can have at most one inflection point, because of the properties of the admissible curves evolving under the curve shortening flow. However, it can not have exactly one inflection point, since otherwise it will fail to be contained inside the strip. Thus, $\chi_t$ is a non-compact convex ancient solution. The only such curve is the Grim Reaper (see \cite{bourni20} for the complete classification of the ancient convex curve shortening flows). $\chi_t$ must be the standard Grim Reaper of width $\pi$, because of the exponential proximity to the walls of Yin-Yang. However, this causes a contradiction, because at the time $t=0$ the limit curve is the standard Grim Reaper while the curves $\beta_0^j$ are no more close enough to Grim Reaper near the tips.

\

{\bf{ \bf{Case 1$^{\circ}$b:} }} When $0 \leq K < \infty$, again, the curves $\beta_t^j$ with $t \in [-K_j,0]$ converge to some limit flow $\chi_t$ with $t \in [-K,0]$. The curves $\beta_{-K_j}^j$, pointed at the tips, converge to the standard Grim Reaper, because they are rescaled admissible initial conditions and, as $j \rightarrow \infty$, near the tips they are congruent to the ever growing chunks of the standard Grim Reaper. So, $\chi_{-K}$ is the standard Grim Reaper and thus $\chi_t$ are also the standard Grim Reaper. Similar to the previous case, this causes another contradiction at the time $t=0$.

\

The above argumentation by contradiction shows that the stability property can not fail for infinitely many indices $j$. {Case 1$^{\circ}$ is studied completely.}

\

{\bf{ \bf{Case 2$^{\circ}$:} }}Now, it will be proved that, when the stability argument holds true, then the proximity argument holds true as well. This will conclude {that the proximity property as well can not fail for infinitely many indices $j$.}

\

First, it will be proved that $\Gamma(s,t)$ is a hyperbolic supersolution.

\

Define the functions $\psi$ and $\eta$:
\begin{equation}
\psi(s,t) := \dfrac{1}{a(t)} \exp[-a(t)(L(t)+s)], \indent \eta(s,t) := \dfrac{1}{a(t)} \exp[-a(t)(R(t)-s)]
\end{equation}
Clearly, $\phi = \psi + \eta$. After computing for $A(\psi)/\psi$ and $A(\eta)/\eta$ it follows that
\begin{eqnarray} \label{Af/f}
\dfrac{A(\psi)}{\psi} = -\dfrac{\dot{a}}{a} - \dot{a}(L+s) -a\dot{L} - a^2 - a \langle F(s),N(s) \rangle - H(s)^2 \\ 
\dfrac{A(\eta)}{\eta} = -\dfrac{\dot{a}}{a} - \dot{a}(R-s) -a\dot{R} - a^2 + a \langle F(s),N(s) \rangle - H(s)^2
\end{eqnarray}

\

Also {due to Lemma \ref{Asymptotics} it can be} assumed that for $t \in [0,T_j)$ (for the tips far enough away from the origin) the mutual ratios of the following quantities are between $1-\epsilon$ and $1+\epsilon$:
\begin{equation} \label{parameter asymptotes}
(3 L_{\pm}(t))^{1/3}, \indent \dfrac{1}{|H(\pm L_{\pm}(t))|}, \indent |F(\pm L_{\pm}(t))|, \indent \dfrac{\pi}{w(- L_{-}(t)))} = \dfrac{\pi}{w(t)} = \dfrac{a(t)}{1-\epsilon}
\end{equation}

\

For the change in time of $L_-(t)$ one has that, $-\dot{L}_-(t) $ is at least the amount of the drift of Yin-Yang due to its rotation (the drift is $\langle -F_t ,T \rangle =  \langle F,N \rangle$, the magnitude of the projection of the velocity to the tangent space) plus the amount of the drift of the tip along the walls, when projected, due to its flow inside the corridor. 

\

The tip is moving in the direction of $-iF(-L_-(t))$ (the direction perpendicular to the radial tangent at the tip) and its speed is at least $\pi(1-\epsilon)/w(t) = a(t)$, because of the approximate Grim Reaper near the tip. The amount of the drift of the tip along the walls is equal to the magnitude of the projection of its velocity to the tangent vector of Yin-Yang.

\

Accounting {the drift term of Yin-Yang and the drift of the tip along the walls and} the asymptotic behaviors of the related quantities in the equation (\ref{parameter asymptotes}) one obtains the following sequence of inequalities:

\begin{eqnarray} \nonumber
-\dot{L}_-(t) &\geq & \langle F(-L_-(t)) , N(-L_-(t)) \rangle + a(t) \dfrac{\langle -iF(-L_-(t)) , T(-L_-(t)) \rangle}{|F(-L_-(t))|} = \\ \label{dL0/dt} 
&=& \langle F,N \rangle + a \dfrac{\langle F,N \rangle}{|F|} = \langle F,N \rangle + a \dfrac{\sqrt{|F|^2 - \langle F,T \rangle ^2}}{|F|} = \langle F,N \rangle + a \sqrt{1 - \dfrac{H^2}{|F|^2}} \geq \\ \nonumber 
&\geq & \langle F,N \rangle + a(1 - (1+\epsilon)H^4) \geq a + \langle F,N \rangle - (1+\epsilon) H^3 
\end{eqnarray}

\

Due to the asymptotics of $w(t) \cong \pi (3L_-(t))^{-1/3}$ and from the equation (\ref{dL0/dt}) it follows that
\begin{eqnarray} \label{dw/dt}
\dot{w}(t) \geq \dfrac{-\pi (1-\epsilon)\dot{L}_-}{(3L_-)^{4/3}} \geq \dfrac{\pi (1-\epsilon)^2}{3L_-} \geq \dfrac{w^3}{\pi^2} = \dfrac{(1-\epsilon)^3}{a^3}
\end{eqnarray}

\

After combining the estimates {in the equations (\ref{dL0/dt}) and (\ref{dw/dt})} and using the asymptotics the following lower bound is obtained:
\begin{eqnarray} \label{adL/dt}
-a \dot{L} \geq a^2 + a \langle F,N \rangle - (1+\epsilon)aH^3 + 100a\dfrac{w^3}{\pi^2} \geq a^2 + a \langle F,N \rangle + \dfrac{10}{a^2}
\end{eqnarray}

\

Similar estimates are true for $-\dot{L}_+(t)$ and $-a\dot{R}(t)$ as well.

\

{By differentiating twice,} it is straightforward to check that the function $a |\langle F(s),N(s) \rangle| + H(s)^2$ is convex and has the absolute minimum at $s=0$. Hence, during the estimates for the lower bounds one may take $s = \pm L_{\pm}(t)$ {for the quantity $a |\langle F(s),N(s) \rangle| + H(s)^2$}.

\

Clearly, because of all the above given relations {(\ref{Af/f}) and (\ref{adL/dt})} one has that 
\begin{eqnarray}
& A(\psi) \geq \dfrac{10}{a^2} \psi - \dot{a}(L+s)\psi \indent \indent \indent A(\eta)  \geq \dfrac{10}{a^2} \eta - \dot{a}(R-s)\eta & \\
& A(\phi) \geq \dfrac{10}{a^2} \phi - \dot{a}(L+s)\psi - \dot{a}(R-s)\eta &
\end{eqnarray}

\

Obviously, $-\dot{a} \geq 0$ and $R-s \geq 0$. Moreover, {because of the estimate in (\ref{dw/dt}) on has that}
\begin{eqnarray}
-\dot{a} = \dfrac{\pi(1-\epsilon)\dot{w}}{w^2} \geq \dfrac{(1-\epsilon)w}{\pi} \geq \dfrac{(1-\epsilon)^2}{a}
\end{eqnarray}



\

Whenever $L+s \gtrsim \frac{\log a}{a}$ or $R-s \gtrsim \frac{\log a}{a}$, the higher order terms in $\phi$ in the equation (\ref{terms}) can be dominated by the linear terms $- \dot{a}(L+s)\psi$  or $- \dot{a}(R-s)\eta$, since $\phi$ is exponentially decaying in the variables $L+s$ or $R-s$ with the decay rate $a$.

\

In the case, when $0 \leq R-s, L+s \lesssim \frac{\log a}{a} \ll 1$ (very close to the endpoints of the barrier curve), effectively, $s=R$ or $s=-L$, $\eta = \phi$ or $\psi = \phi$, $H=\pm\frac{1}{a}$, $\phi_s = \pm a\phi $, $\phi_{ss} = a^2\phi$, and $A(\phi) = \frac{10}{a^2} \phi$.  After plugging the values directly in the inequality (\ref{supersolution condition}), the supersolution condition, and keeping the track of the ``effectively large terms'' one gets that
\begin{eqnarray} \label{higher terms}
A(\phi)+Q(\phi)+P(\phi) \cong \dfrac{10}{a^2}\phi \mp 3a \phi^2 + a^4\phi^3 \mp a^3\phi^4 = \\  \label{long inequality}
= \phi \left( \dfrac{1}{a^2} + \left(\dfrac{3}{a} \mp \dfrac{1}{2}a^2\phi \right) ^2 + a^3\phi^2 \left( \dfrac{3}{4}a \mp \phi \right) \right) \geq 0
\end{eqnarray}

\

It should be noted that, in the equation (\ref{higher terms}) near the endpoints of the $\Gamma$-curve the cubic term $a^4\phi^3$ corresponds to $\phi_s^2 \phi_{ss}$ and is of the order $O(a)$, while the linear term is $O(1/a^3)$ and the other two are $O(1/a)$. This remark justifies the necessity for the careful analysis near the endpoints of the barrier curve and is another way to prove the inequality (\ref{long inequality}). 

\ 

This separate ``endpoint analysis'' concludes that $\Gamma(s,t)$ is indeed a hyperbolic supersolution for $t \in [0,T_j)$.

\

It must further be proved that the curve $\Gamma(s,t)$ lies inside the admissible curve $\alpha_t$ for all times $t \in [0,T_j)$ to completely justify that the $\Gamma$-curve is a barrier curve. This is true, because

\

\begin{enumerate}

\item for the admissible initial condition $\alpha_0$ the exponential decay rate of the glued chunk of Grim Reaper, when viewed as a graph over Yin-Yang in the neighborhood of $tip(\alpha_0)$ until it coincides with the Yin-Yang, is equal to $\pi/w(t=0)$, while for $\Gamma(s,0)$ it is smaller; it is $\pi(1-\epsilon)/w(0)$ to be precise. This will guarantee that $\Gamma(s,0)$ is inside $\alpha_0$ given that $\Vert tip(\alpha_0) \Vert \gg 1$ (the gluing took place far away) and $\lambda_{\pm} > e^{-\lambda}$ (to ensure that the endpoints of $\Gamma(s,0)$ are inside $\alpha_0$);

\

\item by the virtue of the definition of the curve $\Gamma(s,t)$, its endpoints are located in the vicinity of $tip(\alpha_t)$ (the distance $\lambda w(t)$ away from the tip), thus avoiding the intersection with any point away from the vicinity; the endpoints avoid intersecting with any point in the vicinity of the tip, since near the tip the admissible curve $\alpha_t$ looks $\epsilon-$close to Grim Reaper after proper rescaling and therefore the exponential decay rate of $\alpha_t$, when viewed as a graph over Yin-Yang in the neighborhood of the tip, is at least $\pi(1-\epsilon)/w(t) = a(t)$, which is exactly the decay rate of $\Gamma(s,t)$;

\

\item the $\Gamma$-curves are the supersolutions for the curve shortening flow and therefore can not be tangent to the $\alpha$-curves at some point in spacetime.

\end{enumerate}

\

{Case 2$^{\circ}$ is completely analyzed.}

\

{Now the existence of the uniform bound on the curvature function will be proved.} {Assume that for some given $\epsilon>0$ there is no such constant $M_{\epsilon}$. Then there exists a sequence of rescaled admissible curves $\{\alpha_t^j\}_{j=1}^{j=\infty}$ with $t \in [0,t_j]$, where $t_j$ is the first instance when the maximum of the curvature of the curve $\alpha_{t_j}^j$ is at least $j$ (clearly, the maximum is attained away from the neighborhood of its tip, since near the tip the curve looks $\epsilon-$close to the standard Grim Reaper). Similar to the earlier discussion, there are two possibilities: either $\liminf_{j \rightarrow \infty} t_j > 0$ (``long enough flow'') or $\liminf_{j \rightarrow \infty} t_j = 0$ (``rapid increase of curvature''). For any $j$ at time $t_j$ the part of the curve $\alpha_{t_j}^j$ away from its tip is graphical over Yin-Yang, and it is confined between Yin-Yang and the corresponding hyperbolic supersolution. After interpolating between Yin-Yang and the barrier curve one gets uniform in $j$ bound on the $C^1-$norms of the curves $\alpha_{t_j}^j$.}

\

{If $\liminf_{j \rightarrow \infty} t_j > 0$, then, passing to a subsequence if necessary, the uniform in $j$ bound on $C^1-$norm combined with the parabolic regularity will yield a uniform in $j$ bounds on $C^k-$norms. This will contradict with the fact that the maxima of the curvatures of the curves $\alpha_{t_j}^j$ away from the tips is unbounded in $j$.}

\

{In the case of $\liminf_{j \rightarrow \infty} t_j = 0$, the uniform in $j$ bounds of $C^k-$norms of the curves $\alpha_{0}^j$ away from the tips and the uniform in $j$ bound of $C^0-$norm of the curves $\alpha_{t_j}^j$ away from the tips will yield uniform in $j$ bounds of $C^k-$norms of the curves  $\alpha_{\tau}^j$ for $\tau \in [0,t_j]$. This causes another contradiction.} 
Thus, there must exist such constant $M_{\epsilon}$. The theorem is proved. 
\end{proof} 

\section{Appendix}

\subsection{The Asymptotics of $\bar{\Lambda}$}

Consider the circle of radius $R$ centered at the origin and the following Grim Reaper:
\begin{eqnarray}
y = - \log (\cos (x-(R-\pi/2)) ), \indent R-\pi < x < R 
\end{eqnarray}
i.e., Grim Reaper of width $\pi$ with vertical $x=R-\pi$ and $x=R$ asymptotes. Let $0<z<\pi/2$ be such that the point $(R-z,-\log (\cos ( \pi/2 - z ))$ is the intersection point of the circle and the Grim Reaper, which is the closest to the $x$-axis. It can be easily verified that 
\begin{eqnarray} \nonumber
\sqrt{2Rz-z^2} &=& -  \log \sin z \\ \nonumber
\sqrt{2Rz} &\cong & - \log z \\ 
\bar{\Lambda} &\cong & \sqrt{2Rz} \\ \nonumber
e^{\bar{\Lambda}} &\cong & e^{- \log z} = z \\ \nonumber
\bar{\Lambda}^2 e^{\bar{\Lambda}} &\cong & 2R
\end{eqnarray}
when $R \rightarrow \infty$ and hence when $z \rightarrow 0$.

\subsection{Linear and Quadratic Differential Operators $A$ and $Q$}

The computations of the linear differential operator $A$ and the quadratic terms $Q$ are provided here. For simplicity $\Gamma$ is written as $\Gamma = F + \phi N$, where $F = e^{it}F(s)$, $N = e^{it}N(s)$, and $T = F_s$. In the following computations only the linear and quadratic terms in $\phi$ are kept while higher order terms are dropped:
\begin{eqnarray} \nonumber
\Gamma &=& F + \phi N  \\ \nonumber
\Gamma_t &=& F_t + \phi_t N + \phi N_t \\ \nonumber
\Gamma_s &=& (1-\phi H)T + \phi_s N \\ \nonumber
i \Gamma_s &=& (1-\phi H)N-\phi_sT \\ \nonumber
\Gamma_{ss} &=& -(2\phi_sH + \phi H_s)T + (H (1- \phi H) + \phi_{ss})N \\ \nonumber
|\Gamma_s |^2 & = & (1-\phi H)^2 + \phi_s^2 \\
\langle \Gamma_t,i\Gamma_s \rangle &=& 
 \phi_s \langle F,N \rangle + H(1-\phi H) + \phi_t(1-\phi H) + \phi \phi_s \\ \nonumber
\langle F_t,T \rangle &=& \langle iF,T \rangle = -\langle F,N \rangle \\ \nonumber
\langle F_t,N \rangle &=& \langle iF,N \rangle = \langle F,T \rangle = H \\ \nonumber
\langle N_t,T \rangle &=& \langle iN,T \rangle = \langle -T,T \rangle = -1 \\ \nonumber
\langle N_t,N \rangle &=& 0 \\ \nonumber
\langle \Gamma_{ss},i\Gamma_s  \rangle & = & \phi_{ss}(1-\phi H) + \phi_s(2\phi_s H + \phi H_s) + H(1 - \phi H)^2 \\ \nonumber
|\Gamma_s |^2 \langle \Gamma_t,i\Gamma_s \rangle - \langle \Gamma_{ss},i\Gamma_s  \rangle & \cong & \phi_t - \phi_{ss} + \phi_s \langle F(s), N(s) \rangle - \phi H(s)^2 - \\ \nonumber
& - & \phi H [ 3 \phi_t + 2 \phi_s \langle F,N \rangle - 2 \phi H^2 - \phi_{ss} ] + \phi_s [ \phi - \phi_s H - \phi H_s ]
\end{eqnarray}

\

Hence for the linear part it follows that 
\begin{eqnarray}
& A(\phi)  =  \phi_t - \phi_{ss} + \phi_s \langle F(s), N(s) \rangle - \phi H(s)^2 &
\end{eqnarray}

\

Using the relation $H_s = 1 + H \langle F,N \rangle$ and after the careful rearrangements one gets that
\begin{eqnarray}
A(\phi) + Q(\phi) = (1-3\phi H)A(\phi) - H(2 \phi \phi_{ss} + \phi_s^2 + \phi^2 H^2) \label{linear and quadratic}
\end{eqnarray}







\

\

\begin{center}
   \fontsize{10pt}{8pt}\selectfont
    \textsc{Courant Institute of Mathematical Sciences, New York University, 251 Mercer Street, New York, NY 10012}
\end{center}

e-mail address: $jc7518@nyu.edu$

\end{document}